\documentclass[11pt]{amsart}
\usepackage{epsfig,verbatim}
\newtheorem{theorem}{Theorem}
\newtheorem{lemma}[theorem]{Lemma}
\newtheorem{prop}[theorem]{Proposition}
\newtheorem{coro}[theorem]{Corollary}
\newtheorem{conjecture}[theorem]{Conjecture}

\usepackage{amsfonts}
\newcommand{\field}[1]{\mathbb{#1}}
\newcommand{\R}{\field{R}}
\newcommand{\jj}{\mathbf j}
\newcommand{\ii}{\mathbf i}
\newcommand{\hh}{\mathbf h}

\newcommand{\N}{\field{N}}

\newcommand{\reals}{{\R}}

\newcommand{\RR}[1]{\mathbb{R}^{#1}}
\newcommand{\FKG}{(F\,K\!G)}
\newcommand{\GKF}{(G\,K\!F)}

\newcommand{\bfn}{{\bf n}}

\newcommand{\al}{\alpha}

\newcommand{\vt}{\vartheta}
\newcommand{\si}{\sigma}

\newcommand{\eps}{\varepsilon}

\newcommand{\Ksi}{K_{\si,\beta_m,\beta_s}}
\newcommand{\Kwsi}{K_{\beta_m,\beta_s}}
\newcommand{\KH}{K^H_{\beta_m,\beta_s}}
\newcommand{\MM}{\mathcal M}

\newcommand{\hp}{\widehat{p}}
\newcommand{\hq}{\widehat{q}}
\newcommand{\hx}{\widehat{x}}

\newcommand{\abel}[1]{}

\begin{document}
\title{
Asymptotics of certain coagulation-fragmentation
processes and  invariant Poisson-Dirichlet measures}
\author{Eddy Mayer-Wolf,  Ofer Zeitouni, Martin P.\ W.\ Zerner}
\thanks{The work of O.\ Z.\  was
supported in part by the fund for promotion of research at
the Technion,
and in part by a US-Israel
BSF grant. M.\ Z.\ was supported by an Aly Kaufman Fellowship.}
\begin{abstract}
We consider Markov chains on the space of (countable) partitions
of the interval $[0,1]$, obtained first by size biased sampling twice (allowing
repetitions)  and then merging the parts with probability
$\beta_m$ (if the sampled parts are distinct)
or splitting the part with probability $\beta_s$
according to a law $\sigma$ (if the same part was sampled twice).
We characterize invariant probability measures for such chains. In particular,
if $\sigma$ is the uniform measure then the Poisson-Dirichlet law is
an invariant probability measure, and it is unique within a suitably defined
class of ``analytic'' invariant
measures. We also derive transience and recurrence criteria
for these chains.
\end{abstract}
\maketitle
\markboth{EDDY MAYER-WOLF,  OFER ZEITOUNI, MARTIN ZERNER}
{COAGULATION-FRAGMENTATION and POISSON-DIRICHLET}
\section{Introduction and statement of results}

Let $\Omega_1$ denote the space of (ordered) partitions of $1$, that
is
\[\Omega_1:=\left\{p=(p_i)_{i\geq 1}\ :\ p_1\geq p_2\geq
... \geq 0,\ p_1+p_2+\ldots=1\right\}.\]
By {\it size-biased} sampling according to a point $p\in \Omega_1$
we mean picking the $j$-th  part $p_j$ with probability $p_j$.
The starting point for our study is the following Markov chain on $\Omega_1$,
which we call a {\it coagulation-fragmentation} process:
size-bias sample (with replacement) two parts from $p$. If the same
part was picked twice, split it (uniformly), and reorder the partition.
If different parts were picked, merge them, and reorder the partition.

We call this Markov chain the {\it basic chain}.
We first bumped into it
in the context of
triangulation of random Riemann surfaces \cite{bob}. It turns out that it
was already considered in \cite{Tsilevich}, in connection with
``virtual permutations'' and
the Poisson-Dirichlet process. Recall that the Poisson-Dirichlet measure
(with parameter 1) can be described as
 the probability distribution of $(Y_n)_{n\geq 1}$ on
$\Omega_1$ obtained by setting $Y_1=U_1$,
$Y_{n+1}=U_{n+1}(1-\sum_{j=1}^{n} Y_j)$, and reordering the sequence
$(Y_n)_{n}$, where $(U_n)_n$ is a sequence of i.i.d.\
Uniform[0,1] random variables. Tsilevich showed in \cite{Tsilevich} that
the Poisson-Dirichlet distribution is an invariant probability measure for the
Markov chain described above, and raised the question whether
such an
 invariant probability measure is unique. While we do not completely
resolve this question, a corollary of our results (c.f. Theorem
\ref{theo-uniq}) is that the Poisson-Dirichlet law is the unique invariant
measure for the basic chain which satisfies certain regularity conditions.

Of course, the question of invariant probability measure is only one
among many concerning the large time behavior of
the basic chain. Also, it turns out that one may
extend the definition of the basic chain to obtain a
Poisson-Dirichlet measure with any parameter as an
invariant probability measure, generalizing the
result of \cite{Tsilevich}.
We thus consider a slightly more general model, as follows.

For any nonnegative
sequence $x=(x_i)_i$, let $|x|=\sum_i x_i$, the $\ell_1$ norm of $x$, and
$|x|_2^2=\sum_i x_i^2$.
Set
\[\Omega\ =\left\{p=(p_i)_{i\geq 1}\ :\ p_1\geq p_2\geq... \geq 0,\ \ \ \
                                                   \ 0<|p|<\infty\right\}
\]
and $\Omega_{\leq}=\{p\in \Omega:
|p|\leq 1\}$.
Let ${\bf 0} =(0,0,\ldots)$ and define
$\bar \Omega=\Omega\cup \{{\bf 0}\}$ and $\bar \Omega_{\leq}=\Omega_{\leq}\cup
\{{\bf 0}\}$.
Unless otherwise stated, we equip all these spaces
with
the topology
induced from  the product topology on $\R^\N$. In particular,
$\bar \Omega_{\leq}$ is then a compact space.

For a topological space $X$ with
Borel $\sigma$-field  ${\mathcal F}$ we denote by
${\mathcal M}_1(X)$ the set of
all probability measures on $(X,{\mathcal F})$ and equip it with the
topology of weak convergence. ${\mathcal M}_+(X)$ denotes
the space of all (nonnegative) measures on $(X,{\mathcal F})$.

Define the following two operators,
called the {\it merge} and {\it split} operators, on $\bar \Omega$, as follows:
\begin{eqnarray*}
M_{ij}  &:&\bar \Omega\to\bar \Omega,\hspace{.5cm}
    M_{ij}p=\mbox{the nonincreasing sequence obtained by  merging}\\
          & &\hspace*{3.2cm} \ p_i\ \mbox{and}\ p_j\ \mbox{into}\ p_i+p_j\,,
                                                       i\ne j   \\
  S_i^u   &:&\bar \Omega\to\bar \Omega,\hspace{.7cm}
    S_i^u p=\mbox{the nonincreasing sequence obtained by splitting}\ p_i\   \\
          & &\hspace*{3.2cm} \mbox{into}\ \ u\,p_i\ \ \mbox{and}\ \ (1-u)p_i
\,,
  0<u<1
\end{eqnarray*}
Note that the operators $M_{ij}$ and $S_i^u$ preserve the
$\ell_1$ norm.
Let $\si\in{\mathcal M}_1((0,1/2])$ be a probability measure
on $(0,1/2]$ (the {\it splitting measure}).
For $p\in\bar\Omega_\leq$  and $\beta_m,\beta_s\in (0,1]$,
we then consider the Markov process generated in $\bar \Omega_{\leq}$ by the kernel
\begin{eqnarray*}
K_{\si,\beta_m,\beta_s}(p,\cdot)&:=&
2\beta_m\sum_{i<j}p_i p_j
\delta_{M_{ij}p}(\cdot)
+\beta_s\ \sum_i p_i^2 \int \delta_{S_i^u p}(\cdot)
\ d\si(u)\\
&&+\ \left(1-\beta_m |p|^2+
(\beta_m-\beta_s) |p|_2^2
\right) \delta_p(\cdot).
\end{eqnarray*}
It is straightforward to check (see Lemma \ref{cont} below) that
$K_{\sigma,\beta_m,\beta_s}$ is Feller continuous.
The basic chain  corresponds to
$\sigma=U(0,1/2]$, with  $\beta_s=\beta_m=1$.

\begin{figure}[t]
\vspace*{-20mm}
\begin{center}
\psfig{figure=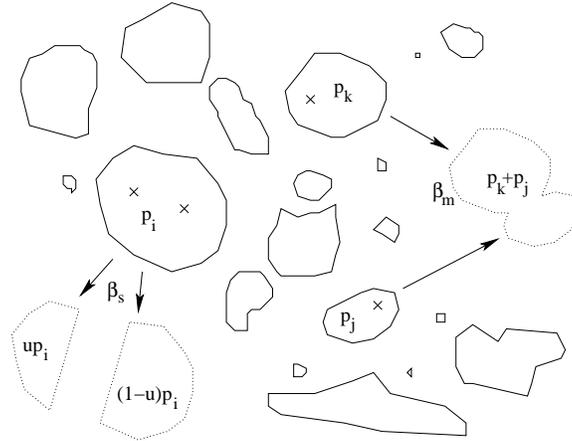,height=280pt}\vspace*{-20mm}
\end{center}
\caption{On the left side a part of size $p_i$ has been chosen twice and
is split with probability
$\beta_s$. On the right side two different parts of sizes $p_i$ and $p_j$ have been
chosen and are merged
with probability $\beta_m$.}
\vspace*{-0mm}
\end{figure}

It is also not hard to check (see Theorem \ref{exists} below) that
there always exists a $K_{\si,\beta_m,\beta_s}$-invariant probability
measure $\mu\in{\mathcal M}_1(\Omega_1)$. Basic properties of any such
invariant probability measure are collected  in Lemma \ref{schteim}
and Proposition \ref{prop}. Our first
result is the following characterization of those kernels that yield
invariant probability measures which are supported on finite (respectively infinite)
partitions. To this end, let $S:=\{p\in\Omega_1\ |\ \exists i\geq 2: p_i=0\}$ be the
set of finite
partitions.
\begin{theorem} {\bf (Support properties)}
\label{theo-duo}
For any $K_{\si,\beta_m,\beta_s}$-invariant
$\mu\in{\mathcal M}_1(\Omega_1)$,
\begin{eqnarray*}
\mu[S]=1&\mbox{\ if\ }&
\int\frac{1}{x}\ d\si(x)\ <\ \infty\qquad\mbox{and}\\
\mu[S]=0&\mbox{\ if\ }&
\int\frac{1}{x}\ d\si(x)\ =\ \infty.
\end{eqnarray*}
\end{theorem}\abel{theo-duo}
Transience and recurrence criteria (which, unfortunately, do
not settle the case $\sigma=U(0,1/2]$!) are provided in the:
\begin{theorem} {\bf (Recurrence and transience)}
\label{rectrans}
\abel{rectrans}
The state
$\bar{p}=(1,0,0,\ldots)$
is  positive recurrent for
$K_{\si,\beta_m,\beta_s}$
if and only if $\int 1/x\ d\si(x)<\infty$.
If however
\begin{equation}\label{bush}
\int_0^{1/2}\frac{1}{\si[(0,x]]}\ dx<\infty
\end{equation}\abel{bush}
then $\bar{p}$ is a transient state for
$K_{\si,\beta_m,\beta_s}$.
\end{theorem}

We now turn to the case $\sigma=U(0,1/2]$.
In order to define invariant probability measures in this case,
set $\pi:\Omega\to\Omega_1,
  \hp:=\pi(p)=({p_i}/{|p|})_{i\ge 1}\,.$
For each $\theta>0$ consider the Poisson process on $\RR{}_+$ with intensity
 measure $\nu_\theta(dx)=\theta x^{-1}{e^{-x}}\,dx$ which can be seen
 either as a Poisson random measure $N(A;\omega)$ on the positive real line or
 as a random variable $X=(X_i)_{i=1}^\infty$ taking values in\ $\Omega$
 whose distribution shall be denoted by $\mu_\theta$, with expectation
 operator $E_\theta$.
 (Indeed, $E_\theta|X|=E_\theta\int_0^\infty x\,N(dx)
     =\int_0^\infty x\,\nu_\theta(dx)<\infty$\
while $P_\theta(|X|=0)=\exp(-\nu_\theta[(0,\infty)])=0$,
     and thus $X\in\Omega$\ a.s.).
 A useful feature of such a Poisson process is that for any Borel subset $A$
 of $\RR{}_+$ with $0<\nu_\theta(A)<\infty$, and conditioned on $\{N(A)=n\}$,
 the $n$ points in $A$ are distributed as $n$ independent variables chosen each
 according to the law $\nu_\theta(\,\cdot\,|A)$.
 The Poisson-Dirichlet measure $\widehat{\mu}_\theta$ on $\Omega_1$ is defined
 to be the distribution of $(\widehat{X}_i)_{i\ge 1}$. In other words,
 $\widehat{\mu}_\theta=\mu_\theta\circ\pi^{-1}$.
In the case $\theta=1$ it coincides with the previously described
Poisson-Dirichlet measure. See~\cite{Ki1},
 \cite{Ki2} and \cite{arratia} for more details and additional properties of Poisson-Dirichlet
processes.

We show in Theorem \ref{theo-uniq} below that, when $\sigma=U(0,1/2]$,
for each choice of  $\beta_m,\beta_s$
 there is a Poisson--Dirichlet measure which is invariant for
 $K_{\sigma,\beta_m,\beta_s}$.
We also show that it is, in this case,
 the unique
invariant probability
 measure in a class $\mathcal A$, which we proceed to define.
Set
\[\Bar \Omega^k_<:=\left\{(x_i)_{1\leq i\leq k}\ :
\ x_i\geq 0,
 x_1+x_2+\ldots+x_k< 1\right\}\]
and denote by $A_k$ the set of real valued functions on $\bar\Omega_<^k$
that coincide (leb$^k$-a.e.)
with a function which has a real analytic extension to some open
neighborhood of $\bar \Omega_<^k$. (Here and throughout,
$\mbox{leb}^k$ denotes the $k$-dimensional Lebesgue measure; all
we shall use is that real analytic functions in a connected domain
can be recovered from their derivatives at an internal point.)
For any $\mu \in {\mathcal M}_1(\Omega_1)$ and
each integer $k$, define the measure
$\mu_k\in{\mathcal M}_+\left(\bar \Omega^k_<\right)$ by
\[\mu_k(B)=E_\mu \left[\sum_{\jj\in\mathbb{N}^k_{\neq}}
     \left(\prod_{i=1}^kp_{j_i}\right)\,{\mathbf 1}_{_B}(p_{j_1},\ldots,p_{j_k})\right]
\,,
       \hspace{1cm}B\in{\mathcal B}_{\Bar \Omega_<^k}\ \,   \]
 (here $\mathbb{N}^k_{\neq}=\left\{\jj\in\mathbb{N}^k\ |\
                            j_i\neq j_{i'}\ \mbox{if}\ i\neq i'\,\right\}$).\
An alternative description of
$\mu_k$ is the following one: pick a random partition $p$ according to $\mu$ and then
sample size-biased
independently (with replacement) $k$ parts $p_{i_1},
\ldots,p_{i_k}$ from $p$.
Then,
$$\mu_k(B)=P(\mbox{the $i_j$-s are pairwise distinct, and}
\, (p_{i_1},
\ldots,p_{i_k})\in B)\,.$$
Part of the proof of part (b) of Theorem~\ref{theo-uniq} below
will consist in verifying that these measures $(\mu_k)_{k\ge 1}$
characterize $\mu$ (see \cite[Th. 4]{pitmanaap} for a similar
argument in a closely related context).

Set for $k\in\N$,
\[ {\mathcal A}_k=\left\{\,\mu\in{\mathcal M}_1(\Omega_1)\,\left|\,
 \mu_k\ll\mbox{leb}^k,
    m_k:= \frac{d\mu_k}{d\,\mbox{leb}^k}\in A_k
\right.\right\}.\]
Our main result is part (b) of the following:
\begin{theorem} {\bf (Poisson-Dirichlet law)}
\label{theo-uniq}
Assume $\sigma=U(0,1/2]$ and
fix $\theta=\beta_s/\beta_m$.\\
(a) The Poisson-Dirichlet law of parameter
$\theta$ belongs to
$\mathcal A:=\bigcap_{k=1}^\infty{\mathcal A}_k$,
and is invariant (in fact: reversing)
for the kernel $K_{\sigma,\beta_m,\beta_s}$.\\
(b)
Assume  a probability measure $\mu\in {\mathcal A}$
is $K_{\si,\beta_m,\beta_s}$-invariant.
Then $\mu$
is the Poisson-Dirichlet law of parameter $\theta$.
\end{theorem}
\abel{theo-uniq}

The structure of the paper is as follows: In Section \ref{sec-pre},
we  prove the Feller property of $K_{\sigma,\beta_m,\beta_s}$, the existence
of invariant probability measures for it, and some of their basic properties.
Section
\ref{sec-support} and \ref{sec-trans} are devoted to the proofs of
Theorems \ref{theo-duo} and  \ref{rectrans} respectively,
Section \ref{sec-Pois} studies the
Poisson-Dirichlet measures and provides the proof of Theorem
\ref{theo-uniq}. We conclude in Section \ref{concluding}
with a list of comments
and open problems.

\section{Preliminaries}
\label{sec-pre}

For fixed $\si\in{\mathcal M}_1((0,1/2])$, $\beta_m,\beta_s\in(0,1]$ and $p\in\bar
\Omega_\leq$  we denote by
$P_p\in{\mathcal M}_1(\bar \Omega_\leq^{\N\cup\{0\}})$ the law of the Markov process on
$\bar
\Omega_\leq$  with kernel
$K_{\si,\beta_m,\beta_s}$ and starting point $p$, i.e.\ $P_p[p(0)=p]=1$.
Whenever $\mu\in{\mathcal M}_1(\bar\Omega_\leq)$,  the law of the
corresponding Markov process with initial
distribution $\mu$ is denoted by $P_\mu$.
 In both cases, we use $(p(n))_{n\geq 0}$ to denote the resulting process.

\begin{lemma}\label{cont}
The kernel $K_{\si,\beta_m,\beta_s}$
is Feller, i.e.\
for
any continuous function
$f:\bar \Omega_\leq\to\R$,
the map $\bar \Omega_\leq\to\R,$
$p\mapsto\int f\ dK_{\si,\beta_m,\beta_s}(p,\cdot)$
is continuous.
\end{lemma}
\abel{cont}
\begin{proof}
 We have
 \begin{eqnarray} \label{Kf}
  \int\,fdK_{\sigma,\beta_m,\beta_s}(p,\cdot)
       &= &2\beta_m\sum_{i=1}^\infty
   p_i\sum_{j=i+1}^\infty p_j\left(f(M_{ij}p)-f(p)\right)\hspace{3cm}
                                              \nonumber\\
            &&\hspace*{1cm}+\beta_s\sum_{i=1}^\infty
      p_i^2\int\left(f(S_i^up)-f(p)\right)d\sigma(u)
            \ \ \ +\ f(p)\nonumber \\
  &=:&2\beta_m\sum_{i=1}^\infty p_ig_i(p)+\beta_s\sum_{i=1}^\infty p_i^2h_i(p)
                                \ +f(p).
 \end{eqnarray}

One may assume that $f(p)$ is of the form $F(p_1,\ldots,p_k)$ with
$k\in\mathbb{N}$
and $F\in C(\bar \Omega_{\le}^k)$, since any
 $f\in C(\bar \Omega_{\le})$
 can be uniformly approximated by such functions,
and denote accordingly $\|p\|_k $
 the $\mathbb{R}^k$ norm of $p$'s first $k$ components.
We shall prove the lemma in
 this case by showing that both sums
in~(\ref{Kf}) contain finitely many nonzero
 terms, this number being uniformly bounded on
some open neighborhood of a given $q$,
 and that $g_i$ and $h_i$ are continuous for every $i$.

 For the second sum these two facts
are trivial: $S_i^up$ and $p$ coincide in their
 first $k$ components  $\forall u\in(0,1/2],\
 \forall i>k$, since splitting a
 component doesn't affect the ordering
of the larger ones, and thus $h_i\equiv 0$
 for $i>k$. Moreover, $h_i$'s continuity
follows from
equicontinuity of  $(S_i^u)_{u\in(0,1)}$.

 As for the first sum, given
$q\in \bar \Omega_{\le}$ with positive components
 (the necessary modification when $q$ has zero components
 is straightforward), let
 $n=n(q)>k$ be such that $q_n<\frac{1}{4}\,q_k$
 and consider $q$'s open neighborhood
 $U=U(q)=\left\{p\in \bar \Omega_{\le}\
                   :\ p_k>\frac{2}{3}q_k,\ p_n<\frac{4}{3}q_n \right\}$.
 In particular, for all $p\in U$,\ \
$p_n<\frac{1}{2}p_k$ and thus, when  $j>i>n$ ,\ \
 $p_i+p_j\le 2p_n<p_k$, which means that $M_{ij}p$
 and $p$ coincide in their first
 $k$ components, or that $g_i(p)=0$ for every $i>n(q)$ and $p\in U(q)$.

 Finally, each $g_i$ is continuous because the series
 defining it converges uniformly.
 Indeed, for $j>i$ and uniformly in $p,\
 \ \ \|M_{ij}p-p\|_k\le p_j\le \frac{1}{j}$.
 For a given $\varepsilon>0$, choose $j_0\in\mathbb{N}$ such that
 $|F(y)-F(x)|<\varepsilon$  whenever $\|y-x\|_k<\frac{1}{j_0}$. Then
 \[ \left|\sum_{j=j_0}^\infty  p_j\left(f(M_{ij}p)-f(p)\right)\right|
              <\varepsilon\sum_{j=j_0}^\infty p_j \le \varepsilon \]
 which proves the uniform convergence.
\end{proof}

\begin{lemma}\label{schteim}
\abel{schteim}
Let $\mu\in{\mathcal M}_1(\bar\Omega_\leq)$
be $K_{\si,\beta_m,\beta_s}$-invariant. Then
\begin{equation}\label{half}
\int |p|_2^2\ d\mu=\frac{\beta_m}{\beta_m+\beta_s}
\int |p|^2\ d\mu.
\end{equation}
\abel{half}
Furthermore,
if we set for $n\geq 1$,
\begin{equation}\label{domi}
\nu_0=\delta_{(1,0,0,\ldots)},\quad
\nu_n=\nu_{n-1} K_{\si,\beta_m,\beta_s}\quad \mbox{and}
\quad \bar{\nu}_n=\frac{1}{n}\sum_{k=0}^{n-1}\nu_k\,,
\end{equation}
then for all $n\geq 1$,
\begin{equation}\label{deu}
\int |p|_2^2\ d\bar{\nu}_n\geq\frac{\beta_m}{\beta_m+\beta_s}
\end{equation}
\end{lemma}
\abel{schteim}\abel{half}\abel{domi}\abel{deu}

\begin{proof}
Let $\eps\in [0,1]$, and
consider  the random variable
\[X_\eps:=\sum_i 1_{\eps<p_i}\]
on $\bar \Omega_\leq$
which counts the intervals longer than $\eps$.
We first prove (\ref{half}).
(The value $\eps=0$ is used in the
subsequent proof of (\ref{deu}).)
Assume that $X_\eps$ is finite which is always the case for $\eps>0$ since on
$\bar \Omega_\leq$,
$X_\eps\leq 1/\eps$ and is also true for
$\eps=0$ if only finitely many $p_i$ are non zero.
Then the expected (conditioned on $p$) increment
$\Delta_\eps$ of
$X_\eps$ after one step of the underlying Markov process is
well-defined. It equals
\begin{eqnarray}\nonumber
\Delta_\eps&=&\beta_m\sum_{i\ne j}
p_i p_j(1_{p_i,p_j\leq \eps<p_i+p_j}-1_{\eps<p_i,p_j})
\nonumber \\
&&+ \ \beta_s \sum_i p_i^2 1_{\eps<p_i}\left(\int 1_{\eps<x p_i}\ d\si(x)-\int
1_{\eps\geq (1-x) p_i}\ d\si(x)\right)\nonumber\\
&=&\beta_m  \sum_{i, j} p_i p_j(1_{p_i,p_j\leq \eps<p_i+p_j}-
1_{\eps<p_i,p_j})\label{flug} \\
&&+ \beta_s \sum_i p_i^2
 1_{\eps<p_i}\left(\si[(\eps/p_i,1/2]]-
\si[[1-\eps/p_i,1/2]]\right)\nonumber\\
&&-
\beta_m \sum_i p_i^2
\left(1_{p_i\leq \eps<2p_i}- 1_{\eps<p_i}\right).\nonumber
\end{eqnarray}
\abel{flug}
The right hand side of (\ref{flug}) converges as $\eps$ tends to 0 to
\begin{equation}
\lim_{\eps\searrow 0}\Delta_\eps
=
-\beta_m|p|^2
+(\beta_m+\beta_s)|p|_2^2.\label{rishon}
\end{equation}
\abel{rishon}
Since $\mu$ is $K_{\si,\beta_m,\beta_s}$-invariant we have
$\int \Delta_\eps\ d\mu=0$ for all $\eps$.
 Now (\ref{half}) follows from (\ref{rishon}) by
dominated convergence since $|\Delta_\eps|\leq 2$.

For the proof of (\ref{deu}) note that for all $n\geq 0$,
$\nu_n$ has  full measure on
sequences $p\in\Omega_1$ for which the number $X_0$ of nonvanishing
components is finite because we start with $X_0=1$ $\nu_0$-a.s.\ and $X_0$
can
increase at most by one in each step.
Given such a  $p\in\Omega_1$, the
expected increment $\Delta_0$ of $X_0$ equals (see (\ref{flug}),
 (\ref{rishon}))
$\Delta_0=-\beta_m+(\beta_m+\beta_s)|p|_2^2$. Therefore for $k\geq 0$,
\[\int X_0\ d\nu_{k+1}-\int X_0\ d \nu_k=-
\beta_m +(\beta_m+\beta_s)\int |p|_2^2\ d\nu_k.\]
Summing over $k=0,\ldots,n-1$ yields
\begin{equation}\label{liba}
\int X_0\ d\nu_n-\int X_0\ d\nu_0=
-n\beta_m +(\beta_m+\beta_s)\sum_{k=0}^{n-1}\int |p|_2^2\ d\nu_k.
\end{equation}\abel{liba}
The left hand side of (\ref{liba}) is nonnegative
due to $\int X_0\ d\nu_0=1$ and
$\int X_0\ d\nu_n\geq 1$. This proves (\ref{deu}).
\end{proof}

\begin{theorem}  \label{exists}
There exists a $K_{\si,\beta_m,\beta_s}$-invariant probability
measure $\mu\in{\mathcal M}_1(\Omega_1)$.
\end{theorem}
\abel{exists}
\begin{proof}
Define $\nu_n$ and $\bar{\nu}_n$ as in (\ref{domi}).
Since $\bar \Omega_\leq$ is
compact,
${\mathcal M}_1(\bar \Omega_\leq)$
is compact. Consequently, there are
 $\mu\in{\mathcal
M}_1(\bar \Omega_\leq)$ and a strictly increasing
sequence $(m_n)_n$ of positive integers
such that $\bar{\nu}_{m_n}$ converges weakly towards $\mu$ as $n\to\infty$.
This limiting measure $\mu$ is invariant under $K_{\si,\beta_m,\beta_s}$ by the
following standard argument. For any
continuous function $f:\bar \Omega_\leq\to\R$,
\begin{eqnarray*}
\int f\ d(\mu K_{\si,\beta_m,\beta_s})
&=&
\int\int f\ dK_{\si,\beta_m,\beta_s}(p,\cdot)\ d\mu(p)\\
&=&\lim_{n\to\infty}\int\int f\ dK_{\si,\beta_m,\beta_s}(p,\cdot)\
 d\bar{\nu}_{m_n}(p)\qquad[\mbox{Lemma
\ref{cont}}]\\
&=&\lim_{n\to\infty}\frac{1}{m_n}\sum_{k=0}^{m_n-1}
\int\int f\ dK_{\si,\beta_m,\beta_s}(p,\cdot)\ d\nu_k(p)\\
&=&\lim_{n\to\infty}\frac{1}{m_n}\sum_{k=0}^{m_n-1}
\int f(p)\  d\nu_{k+1}(p)
=\lim_{n\to\infty} \int f d\bar{\nu}_{m_n}=\int f d\mu\,.
\end{eqnarray*}
Hence it remains to show that $\Omega_1$ has full $\mu$-measure,
i.e.\ $\mu[|p|=1]=1$.
To prove this  observe that $|p|_2^2$
(unlike $|p|$) is a continuous function on
$\bar\Omega_\leq$.
Therefore by (\ref{half}), weak convergence and (\ref{deu}),
\[1\geq \int  |p|^2\ d\mu =\frac{\beta_m+\beta_s}{\beta_m} \int |p|_2^2\ d\mu
=\frac{\beta_m+\beta_s}{\beta_m} \lim_{n\to\infty} \int |p|_2^2\ d\bar{\nu}_{m_n}\geq
1\]
by which the first inequality is an equality, and thus $|p|=1\ \mu-a.s.$
\end{proof}

\begin{prop}
\label{prop}
\abel{prop}
If $\mu\in{\mathcal M}_1(\Omega_1)$ is
$K_{\si_i,\beta_{m,i},\beta_{s,i}}$-invariant
for $i=1,2$, then $\si_1=\si_2$ and
$\theta_1:=\beta_{s,1}/\beta_{m,1}=\beta_{s,2}/\beta_{m,2}=:\theta_2$.
\end{prop}
\begin{proof}
Let $k\geq 1$ be an integer
and $\alpha\in\{1,2\}$.
Given $p$, consider the expected increment
 $\Delta_{\alpha,k}$ of $\sum_i p_i^k$
after one step of the process driven by
$K_{\si_\al,\beta_{m,\alpha}\beta_{s,\alpha}}$:
\begin{eqnarray*}
\Delta_{\alpha,k}&=&\beta_{m,\alpha}
\sum_{i\ne j}p_i p_j \left(-p_i^k-p_j^k+(p_i+p_j)^k\right)\\
&&+\beta_{s,\alpha}
\ \sum_i p_i^2\left(-p_i^k+\int(t p_i)^k+((1-t)p_i)^k\ d\si_\alpha(t)\right)
\,.
\end{eqnarray*}
Note that $\int \sum_i p_i^k\ d\mu$ is finite because of
$k \geq 1$. Therefore,
 by invariance,
$\int \Delta_{\alpha,k}\ d\mu=0$, which implies
\[
\beta_{s,\alpha}\left[\int (t^k+(1-t)^k)\ d\si_\alpha(t)-1\right]=\frac{
\beta_{m,\alpha}\int
\sum_{i\ne j}p_i p_j \left(p_i^k+p_j^k-(p_i+p_j)^k\right)\ d\mu}
{\int \sum_ip_i^{2+k}\ d\mu}.
\]
Hence, for any $k$,
\[\frac{\int (t^k+(1-t)^k)\ d\si_1(t)-1}
{\int (t^k+(1-t)^k)\ d\si_2(t)-1}=
\frac{\beta_{m,1}\beta_{s,2}}{\beta_{m,2}\beta_{s,1}}=:\gamma\,.\]
Taking $k\to\infty$ we conclude that $\gamma=1$. This proves the second claim.
In addition, we have
\begin{equation}\label{shalom}
\int (t^k+(1-t)^k)\ d\si_1(t)=\int (t^k+(1-t)^k)\ d\si_2(t)
\end{equation}\abel{shalom}
for all $k\geq 1$. Obviously, (\ref{shalom}) also holds true for $k=0$.
Extend $\sigma_\alpha$ to probability
measures on $[0,1]$ which are
supported on $[0,1/2]$. It is enough for the proof of
$\si_1=\si_2$  to show that
 for all continuous real valued
functions $f$ on $[0,1]$ which vanish on $[1/2,1]$
the integrals $\int f(t)\ d\si_\al(t)$ coincide for  $\al=1,2$.
Fix such an $f$ and choose
a sequence of polynomials
\[\pi_n(t)=\sum_{k=0}^n c_{k,n}t^k\qquad(c_{k,n}\in\R)\]
which converges uniformly on $[0,1]$ to $f$ as $n\to\infty$.
Then $\pi_n(t)+\pi_n(1-t)$ converges
uniformly on $[0,1]$ to $f(t)+f(1-t)$.
Since $f(1-t)$ vanishes on the support of
$\si_1$ and $\si_2$ we get for $\alpha=1,2$,
\begin{eqnarray*}
\int f(t)\ d\si_\alpha(t)&=& \int f(t)\ d\si_\alpha(t)+\int f(1-t)\
d\si_\alpha(t)\\
&=&\lim_{n\to\infty} \sum_{k=0}^n c_{k,n} \int (t^k+(1-t)^k)\ d\si_\alpha(t)
\end{eqnarray*}
which is the same for $\alpha=1$ and $\alpha=2$ due to (\ref{shalom}).
\end{proof}

\newpage
\section{Support properties}
\label{sec-support}

Theorem \ref{theo-duo} is a consequence of  the following result.
\begin{theorem}\label{bit}
Let $\mu\in{\mathcal M_1}(\Omega_1)$ be
$K_{\sigma,\beta_m,\beta_s}$-invariant
and denote $\bar{p}:=(1,0,0,\ldots)$ and
     $(p(n))$'s stopping time
     $H:=\min\{n\ge 1\ :\ p(n)=p(0)\}$.
Then
\[\int\frac{1}{x}\ d\si(x)<\infty\Longleftrightarrow
\mu[S]=1 \Longleftrightarrow
\mu[S]>0 \Longleftrightarrow \mu[\{\bar{p}\}]>0
                             \Longleftrightarrow E_{\bar{p}}[H]<\infty.
\]
\end{theorem}\abel{bit}

\begin{proof}
We start by proving that  $\int1/x\ d\si(x)<\infty $ implies $\mu[S]=1$.
Fix an arbitrary $0<\vt\leq 1/2$ and
 consider the random variables
\[W_n:=\sum_{i\geq 1}p_i 1_{\vt^n<p_i}\quad(m\geq 1).\]
After one step of the process $W_n$ may increase, decrease or stay unchanged.
If we merge two intervals then $W_n$ cannot decrease, but may increase by
the mass
 of one or two intervals
which are smaller than $\vt^n$ but
 become part of an interval which is bigger than
$\vt^n$. If we split an interval then $W_n$ cannot
increase, but it decreases if the original interval
was larger than $\vt^n$ and at least one
of its parts is smaller than $\vt^n$.
Thus given $p$, the expected increment $\Delta$ of $W_n$
after one step of the process is

\begin{eqnarray*}
\Delta
&:=& \Delta_+ - \Delta_-,\qquad \mbox{where}\\
\Delta_+&:=& \beta_m \sum_{i\ne j}p_i p_j
\left(p_i1_{p_i\leq \vt^n<p_j}+p_j1_{p_j\leq \vt^n<p_i}+
(p_i+p_j)1_{p_i,p_j\leq\vt^n<p_i+p_j}\right)\qquad\mbox{and}\\
\Delta_-
&:=& \beta_s \sum_{i}p_i^2 \int \left(p_i 1_{(1-x)p_i\leq \vt^n<p_i}+x p_i
1_{x p_i\leq\vt^n<(1-x)p_i}\right)\ d\si(x).
\end{eqnarray*}
We bound $\Delta_+$ from below by
\begin{eqnarray*}
\Delta_+&\geq& 2\beta_m
\sum_{i, j}p_i^2p_j 1_{p_i\leq \vt^n}\cdot 1_{\vt^n<p_j}\\
&\geq&2\beta_m
\left(\sum_i p_i^2 1_{\vt^{n+1}<p_i\leq\vt^n}\right)\left(\sum_jp_j
1_{\vt^n<p_j}\right)\\
&\geq&2\beta_m
\vt^{2n+2} W_n \sharp I_{n+1}
\end{eqnarray*}
where
\[I_n:=\left\{i\geq 1\ :\ \vt^n<p_i\leq \vt^{n-1}\right\}\qquad(m\geq 1)\]
and $\Delta_-$ from above by
\begin{eqnarray*}
\Delta_-&\leq&
\beta_s \sum_{i\geq 1}\int \left(p_i^3 1_{\vt^n<p_i\leq \vt^n/(1-x)}
+p_i^3 x 1_{\vt^n<p_i}\cdot 1_{p_i\leq \vt^n/x}\right)\ d\si(x)\\
&\leq&\beta_s \sum_{i\geq 1}p_i^3 1_{\vt^n<p_i\leq \vt^{n-1}}\qquad
\mbox{[since $\vt\leq 1/2\leq 1-x$]}\\
&&+\ \beta_s
\sum_{i\geq 1}\int \sum_{j=0}^{n-1}p_i^3 x 1_{\vt^{n-j}<p_i\leq \vt^{n-j-1}}
1_{p_i\leq \vt^n/x}\ d\si(x)\\
&\leq& \beta_s \sum_{i\geq 1}\vt^{3(n-1)} 1_{\vt^n<p_i\leq \vt^{n-1}}\\
&&+\beta_s\ \sum_{i\geq 1}\int
\sum_{j=0}^{n-1}\vt^{3(n-j-1)} x 1_{\vt^{n-j}<p_i\leq \vt^{n-j-1}}
1_{x\leq \vt^j}\ d\si(x)\\
&\leq&\beta_s \vt^{3(n-1)} \sharp I_n
+\beta_s \
\sum_{j=0}^{n-1}\sum_{i\geq 1}\vt^{3(n-j-1)}
 \vt^j 1_{\vt^{n-j}<p_i\leq \vt^{n-j-1}}
\si[(0,\vt^j]]\\
&\leq&\beta_s \vt^{3(n-1)} \sharp I_n
+\beta_s \vt^{3(n-1)}\sum_{j=0}^{n-1}  \vt^{-2j}\si[(0,\vt^j]]\sharp I_{n-j}\\
&\leq&2
\beta_s
\vt^{3(n-1)}\sum_{j=0}^{n-1}\vt^{-2j}\si[(0,\vt^j]] \sharp I_{n-j}.
\end{eqnarray*}
Since $\mu$ is invariant by assumption,
$0=\int\Delta\ d\mu=\int \Delta_+\ d\mu-\int\Delta_-\ d\mu$ and
therefore
\begin{eqnarray*}
2\beta_m
\int W_n \sharp I_{n+1}\ d\mu &\leq& 2\beta_s \vt^{3n-3-2n-2}\sum_{j=0}^{n-1}
\vt^{-2j}\si[(0,\vt^j]]
 \int \sharp I_{n-j}\ d\mu\\
& =& 2 \beta_s
\vt^{-5}\sum_{j=0}^{n-1}\vt^{-j}\si[(0,\vt^j]]\int \vt^{n-j}\sharp I_{n-j}\
d\mu.
\end{eqnarray*}
Consequently\,,
\begin{eqnarray*}
\lefteqn{
\sum_{n\geq 1}\mu[W_n \sharp I_{n+1}\geq 1/2]}\\
&\leq&
\sum_{n\geq 1}2\int W_n \sharp I_{n+1}\ d\mu
\ \leq\ \frac{2\vt^{-5}\beta_s}{\beta_m}
\sum_{n\geq 1}\sum_{j=0}^{n-1}\vt^{-j}\si[(0,\vt^j]]
\int \vt^{n-j}\sharp I_{n-j}\ d\mu\\
&=&\frac{2 \vt^{-5}\beta_s}{\beta_m}
\left(\sum_{j=0}^\infty\vt^{-j}\si[(0,\vt^j]]\right) \sum_{n\geq 1}
\int \vt^{n}\sharp I_{n} \ d\mu\\
 &\leq&
\frac{2\vt^{-5}\beta_s}{(1-\vt)\beta_m}
\left(\sum_{j=0}^\infty(\vt^{-j}-\vt^{-j+1})
\si[(0,\vt^j]]\right)
\sum_{n\geq 1} \int \sum_{i\in I_n}p_i\ d\mu\\
&=&\frac{2\vt^{-5}\beta_s}{(1-\vt)\beta_m}
\left(\int\sum_{j=0}^\infty 1_{\vt^{-j}\leq 1/x}
(\vt^{-j}-\vt^{-j+1})\ d\si(x)\right) \int |p|\ d\mu\\
& \leq&
\frac{2\vt^{-5}\beta_s}{(1-\vt)\beta_m}\int\frac{1}{x}\ d\si(x)
\end{eqnarray*}
which is finite by assumption. Thus by Borel-Cantelli,
 $W_n \sharp I_{n+1} $ is $\mu$-a.s.\ eventually
(for large $n$) less
than 1/2. However, $W_n$ converges $\mu$-a.s.\  to 1 as $n$ tends to $\infty$.
Thus even $\sharp I_{n+1}$ is  $\mu$-a.s.\ eventually less
than 1/2, which means that $I_{n+1}$ is $\mu$-a.s.\ eventually empty, that is
$\mu[S]=1$.

Now we assume $\mu[S]>0$ in which case there exist some $i\geq 1$ and $\eps>0$
      such that $\delta:=\mu[p_i>\eps, p_{i+1}=0]>0.$
      By $i$ successive merges of the positive parts and $\mu$'s invariance
      we obtain
     \begin{equation}\label{bump}
      \mu[\{\bar{p}\}]=\mu[p_1=1]\geq (2\beta_m \eps^2)^{i-1}\delta >0.
     \end{equation}

  Next, we assume $\mu[\{\bar{p}\}]>0$ and note that $K_{\si,\beta_m,\beta_s} 1_S
     = 1_S$ and thus, defining $\bar \mu:=\mu/\mu[S]$,
     one obtains an invariant measure
     supported on $S$. The chain determined by $K_{\si,\beta_m,\beta_s}$ on $S$
     is $\delta_{\bar{p}}$-irreducible,
and has $\bar \mu$
     as invariant measure, with $\bar \mu[\{\bar{p}\}]>0$.
Therefore,
Kac's recurrence theorem~\cite[Theorem 10.2.2]{Meyn} yields $E_{\bar{p}}[H]<\infty$.

Finally, we assume $E_{\bar{p}}[H]<\infty$ and show $\int 1/x\ d\si(x)<\infty$.
If $A:=\{\bar{p}=p(0)\ne p(1)\}$, then $P_{\bar{p}}[A]=\beta_s>0$, and when
    $p\in A$ we write $p(1)=p^{\xi}:=(1-\xi,\xi,0,\ldots)$,
    where $\xi$ has distribution $\sigma$.
Furthermore, restricted to $A$ and conditioned on $\xi$,\ \  $H\ge\tau\ \
 P_{p^\xi}$--a.s., where in terms of the chain's sampling and merge/split
 interpretation, $\tau$ is the first time a marked part of size $\xi$ is sampled,
 i.e. a geometric random variable with parameter $1-(1-\xi)^2\le 2\xi$. Thus
 \[ \infty>E_{\overline{p}}[H]\ge P_{\overline{p}}[A]E_{\overline{p}}[H|A]
      \ge \beta_s\left(1+\int E_{p^\xi}[\tau]d\sigma(\xi)\right)
      \ge \beta_s\left(1+\int \frac{1}{2\xi}d\sigma(\xi)\right). \]

\end{proof}
\begin{coro}
\label{coro}
If $\int1/x\ d\si(x)<\infty$ then there exists a unique
$K_{\si,\beta_m,\beta_s}$-invariant probability measure
$\mu\in{\mathcal M}_1(\Omega_1)$.
\end{coro}
\begin{proof}
In view of Theorem \ref{theo-duo}, for the study of invariant measures
it is enough to restrict attention to the state space
$S$, where the Markov chain $(p(n))_n$ is
$\delta_{\bar{p}}$-irreducible, implying, see
\cite[Chapter  10]{Meyn}, the uniqueness of the invariant measure.
\end{proof}

\section{Transience and recurrence}
\label{sec-trans}
\begin{proof}[Proof of Theorem \ref{rectrans}]
The statement about positive recurrence is included in Theorem \ref{bit}.

The idea for the proof of the transience statement is to show that
under (\ref{bush}) the event that the size of the smallest positive part of the partition never
increases has positive probability. By
\[n_0:=0\quad\mbox{and}\quad n_{j+1}:=\inf\{n>n_j:\ p(n)\ne p(n-1)\}\quad(j\geq 0)\]
we enumerate  the times $n_j$ at which the value of the Markov chain
changes.
Denote by $s_n$ the (random) number of instants among the first
$n$ steps of the Markov chain in which some interval is split.
Since $j- s_{n_j}$
is
the number of steps among the first $n_j$ steps in which two parts are merged and since
this number can never exceed
$s_{n_j}$ if $p(0)=\bar{p}$,
we have that $P_{\bar{p}}$-a.s.,
\begin{equation}\label{as}
s_{n_j}\geq \left\lceil\frac{j}{2}\right\rceil\quad \mbox{for all}\quad j\geq 0.
\end{equation}\abel{as}
Let $(\tau_l)_{l\geq 1}$ denote the times at which some part is split.
This part is split
 into two parts of sizes $\ell(l)$ and $L(l)$ with
$0<\ell(l)\leq L(l)$. According to the model the random variables
$\xi_l:=\ell(l)/(\ell(l)+L(l)),\
l\geq 1,$ are independent with common distribution $\si$.
Further, for any deterministic sequence
$\xi=(\xi_n)_n$, let $P_{\xi,\bar p}[\ \cdot\ ]$
denote the law of the process which evolves using the kernel
$K_{\sigma,\beta_m,\beta_s}$ except that at the times
$\tau_l$ it uses the values $\xi_l$  as the splitting variables. Note that
$$P_{\bar p}[\ \cdot\ ]=\int P_{\xi,\bar p}[\ \cdot\ ]\ d\sigma^{\N}(\xi)\,.$$

Now denote by $q(n):=\min\{p_i(n):\ i\geq 1,\ p_i(n)>0\}\quad (n\geq 0)$
 the size of the smallest positive
 part at
time $n$. We prove that for  $N\geq 0$,
\begin{equation}\label{impl}
q(0)\geq\ldots\geq q(N)\quad\mbox{implies}\quad q(N)\leq
\xi_1\wedge\xi_2\wedge\ldots\wedge\xi_{s_N}.
\end{equation}\abel{impl}
(Here and in the sequel, we take $\xi_1\wedge\ldots\wedge\xi_{s_N}=\infty$
if $s_N=0$). Indeed, we need only consider the case $s_N>0$, in which case
there exists a $1\leq t\leq s_N$
such that $\xi_t=\xi_1\wedge\ldots\wedge\xi_{s_N}$,
and $\tau_t\leq N$. But clearly $q(\tau_t)\leq \xi_t$, and
then the condition $q(1)\geq \cdots\geq q(N)$ and
the fact that $\tau_t\leq N$ imply
$q(N)\leq q(\tau_t)\leq \xi_t=\xi_1\wedge\ldots\wedge\xi_{s_N}$,
as claimed.

Next, fix some $\eps\in
(0,\beta_0/2]$ where $\beta_0:=\min\{\beta_m, \beta_s\}$.
We will prove by induction over $j\geq 1$ that
\begin{equation}\label{HNPP}
P_{\xi,\bar{p}}[\eps>q(1), q(0)\geq \ldots\geq q(n_j)]\geq
\beta_s 1_{\xi_1<\eps}\prod_{k=1}^{j-1}\left(1-\frac{\xi_1\wedge\ldots\wedge
\xi_{\lceil k/2\rceil}}{\beta_0}\right).
\end{equation}\abel{HNPP}
For $j=1$  the left hand side of (\ref{HNPP}) equals the probability that the
unit interval is split
 in the first step with the smaller part being smaller than
$\eps$ which equals $\beta_s 1_{\xi_1<\eps}$.  Assume that (\ref{HNPP}) has
been proved
up to $j$. Then, with $ {\mathcal F}_{n_j}=\sigma(p(n),n\leq n_j)$,
\begin{eqnarray}\nonumber\lefteqn{
P_{\xi,\bar{p}}[\eps>q(1),\ q(0)\geq \ldots\geq q(n_{j+1})]}\\
&=&E_{\xi,\bar{p}}\left[
P_{\xi,\bar{p}}[q(n_j)\geq
q(n_{j+1})\ |\ {\mathcal F}_{n_j}],
\eps>q(1),\ q(0)\geq \ldots\geq q(n_j)\right].\label{wer}
\end{eqnarray}\abel{wer}
Now choose $k$ minimal such that $p_k(n_j)=q(n_j)$. One possibility to achieve
$q(n_j)\geq
q(n_{j+1})$ is not to merge the  part $p_k(n_j)$ in the next step in
which the Markov chain moves.
The probability to do this is
\begin{eqnarray*}
1-\frac{\beta_m\sum_{a:a\ne k}p_a(n_j)p_k(n_j)}{\beta_m
\sum_{a\ne b}p_a(n_j)p_b(n_j)
+\beta_s\sum_a p_a^2(n_j)}
&\geq& 1-\frac{\beta_m q(n_j)
\sum_{a}p_a(n_j)}{\beta_0\sum_{a,b}p_a(n_j)p_b(n_j)}\\
&\geq&1-\frac{q(n_j)}
{\beta_0}.
\end{eqnarray*}
Therefore (\ref{wer}) is greater than or equal to
\[
 E_{\xi,\bar{p}}\left[(1-q(n_j)/\beta_0), \eps>q(1),\ q(0)\geq \ldots\geq q(n_j)
\right].
\]
By (\ref{impl}) this can be estimated from below by
\[
 E_{\xi,\bar{p}}\left[
(1-(\xi_1\wedge\ldots\wedge \xi_{s_{n_j}})/\beta_0),\
\eps>q(1),\ q(0)\geq \ldots\geq q(n_j)
\right].
\]
This is due to (\ref{as}) greater than or equal to
\[(1-(\xi_1\wedge\ldots\wedge \xi_{\lceil j/2 \rceil})/\beta_0)\
 P_{\xi,\bar{p}}\left[\eps>q(1),\ q(0)\geq \ldots\geq q(n_j)\right].
\]
Along with the induction hypothesis this implies (\ref{HNPP}) for $j+1$.\vspace*{1mm}

Taking expectations with respect to $\xi$ in (\ref{HNPP}) yields
\begin{equation}\label{warum}
P_{\bar{p}}[q(n)\leq \eps \quad\mbox{for all $n\geq 1$}]
\geq
E_{\bar{p}}\left[\beta_s1_{\xi_1<\eps}\prod_{k\geq
1}\left(1-\frac{\eps\wedge\xi_2\wedge\ldots\wedge
\xi_{\lceil k/2\rceil}}{\beta_0}\right)\right].
\end{equation}\abel{warum}
By independence of $\xi_1$ from  $\xi_i,\ i\geq 2,$ the right hand side of
(\ref{warum}) equals
\begin{equation}\label{stoned}
\beta_s \left(1-\frac{\eps}{\beta_0}\right)^2P[\xi_1<\eps]
E_{\bar{p}}\left[\prod_{k\geq
2}\left(1-\frac{\eps\wedge\xi_2\wedge\ldots\wedge
\xi_k}{\beta_0}\right)^2\right].
\end{equation}\abel{stoned}
Observe that (\ref{bush}) implies $P[\xi_1<\eps]=\si[(0,\eps)]>0$.
By Jensen's inequality and monotone convergence,
(\ref{stoned}) can be estimated from below by
\[c_1\exp\left(\sum_{k\geq 2} 2 E_{\bar{p}}\left[\ln\left(1-
\frac{\eps\wedge\xi_2\wedge\ldots\wedge
\xi_k}{\beta_0}\right)
\right]\right)
\]
with some positive constant $c_1=c_1(\eps)$. Since $\ln(1-x)\geq -2 x$ for $x\in[0,1/2]$ this is greater than
\begin{equation}\label{oje}
c_1\exp\left(-\frac{4}{\beta_0} \sum_{k\geq 2}E_{\bar{p}}[\xi_2\wedge\ldots\wedge
\xi_k]\right)=c_1\exp\left(-\frac{4}{\beta_0} \int_0^{1/2}
\frac{P_{\bar{p}}[\xi_1> t]}{P_{\bar{p}}[\xi_1\leq t]}\ dt\right)
\end{equation}\abel{oje}
where we used that due to independence
\[E_{\bar{p}}[\xi_2\wedge\ldots\wedge
\xi_k]=\int_0^{1/2}P_{\bar{p}}[\xi_1> t]^{k-1}\ dt.\]
Due to assumption (\ref{bush}), (\ref{oje}) and therefore also the left hand side of
(\ref{warum}) are positive. This implies transience of $\bar{p}$.
\end{proof}

\section{Poisson-Dirichlet invariant probability measures}
\label{sec-Pois}
Throughout this section, the splitting measure is the uniform measure
on $(0,1/2]$. To emphasize this, we use $K_{\beta_m,\beta_s}$ instead
of $K_{\sigma,\beta_m,\beta_s}$ throughout.
Recall that $\theta=\beta_s/\beta_m$.

 \noindent
 It will be convenient to equip $\Omega$
 (but not $\bar \Omega_{\leq}$) with the $\ell_1$ topology
(noting that the Borel $\sigma$-algebra is not affected by this change
of topology),
and to replace the kernel $\Kwsi$  by
 \begin{eqnarray}
  \KH(p,\cdot)   \label{KH}
       = \beta_m\,\sum_{i\neq j}\hp_i \hp_j \delta_{M_{ij}p}(\cdot)
          &+&\beta_s\,\sum_i \hp_i^{\,2}
                   \int_0^1\delta_{S_i^u p}(\cdot)\,du \nonumber\\
          &+&\left( 1-\beta_m
               +(\beta_m-\beta_s)\,|\hp|_2^{\,2} \right) \delta_p(\cdot).
 \end{eqnarray}
 Both kernels coincide on $\Omega_1$ (not on $\Omega_{\le}$).
 However, $\KH$ has the advantage that it is well defined on all of $\Omega$ and
 is homogeneous (hence the superscript $H$) in the sense of the first of the
 following two lemmas,
whose proof is straightforward and in which by a slight abuse
 of notation $\KH$ will denote both the kernel in $\Omega_1$ and in $\Omega$ and also
 the operators induced by these kernels, the distinction being clear from the context.
 \begin{lemma} \label{homogeneous}
  For all $p\in\Omega$,\,
  $\KH(\pi p,\cdot)=\KH(p,\cdot)\circ \pi^{-1}$.\\
  More generally, denoting $(\Pi f)(p)=f(\pi(p))
  $,\ \ \ we have\ \ \ $\KH \Pi=\Pi \KH$.\\

  \noindent In particular, if $\mu\in\MM_1(\Omega)$ is invariant (resp.
  reversing) for $\KH$\ \  then\\  $\mu\circ \pi^{-1}\in\MM_1(\Omega_1)$ is
  invariant (resp. reversing) for  $\Kwsi$.
  \end{lemma}
  \begin{lemma} \label{Feller}
   The kernel $\KH$ maps continuous bounded functions to continuous bounded
functions.
  \end{lemma}
  \begin{proof}[Proof of Lemma \ref{Feller}]
Note that we work with the $\ell_1$ topology, and hence have to modify
the proof in Lemma \ref{cont}. The $\ell_1$ topology  makes the mapping
$p\mapsto \hp$ continuous (when $\Omega_1$ is  equipped
with the induced $\ell_1$ topology).
 Fix $F\in C_b(\Omega)$. By~(\ref{KH}) we have
   \begin{eqnarray} \label{KHF}
     \KH\,F(p)
       &=&\beta_m\sum_{i\ne j}\widehat{p}_i\widehat{p}_j F\left(M_{ij}p\right)
              +\beta_s\sum_i\left(\widehat{p}_i\right)^2
                   \int_0^1F\left(S_i^up\right)\,du  \nonumber\\
       & &\hspace*{3.45cm}+\left(1-\beta_m+(\beta_m-\beta_s)
             |  \hp|_2^{\,2}\right)F(p)\nonumber\\
       \hspace*{2cm}&=&\beta_mK_1(p)+\beta_sK_2(p)+K_3(p).\hspace{4cm}
   \end{eqnarray}
   Note that for $l=1,2\ ,\  K_l(p)$ is of the form
   $\langle T_l(p)\hp,\hp\,\rangle$,
   with $T_l(\cdot)\in C\left(\Omega;L(\ell_1,\ell_\infty)\right)$,
   and $\langle\cdot,\cdot\rangle$ denoting the standard duality pairing.
   In stating this we have used the facts that
$ F$ is continuous and
bounded,
   and that all the mappings $M_{ij}$ and $S_i^u$ are
   contractive.

   The continuity of $K_l,\ \ l=1,2$\ ,\ then follows from
   \[ \langle T_l(q)\hq,\hq\,\rangle-\langle T_l(p)\hp,\hp\,\rangle
         =\langle T_l(q)\hq,\hq-\hp\,\rangle
                              +\langle\left(T_l(q)-T_l(p)\right)\hq,\hp\,\rangle
                              +\langle T_l(p)(\hq-\hp),\hp\,\rangle \]
   after observing that
$|\hq|$ and $\|T_l(q)\|$ remain bounded in any $\ell_1$
   neighborhoods of $p$.

The continuity of $K_3$ is obvious being the product of two
continuous functions of $p$. It has thus been shown that $\KH\,F\in C(\Omega)$.
\end{proof}

\noindent
  \begin{theorem} \label{reversible}
   The Poisson-Dirichlet measure $\widehat{\mu}_\theta\in\MM_1(\Omega_1)$
   is reversing for $\Ksi$ with $\si=U(0,1/2]$.
  \end{theorem}
\abel{reversible}
 \begin{proof}
   By Lemma~\ref{homogeneous} it suffices to verify that
   $\mu_\theta\in \MM_1(\Omega)$ is reversing for the kernel $\KH$,
   which for simplicity will be denoted by $K$ for the rest of this proof.

   We thus need to show that
   \begin{equation} \label{GKF=FKG}
    E_\theta\GKF=E_\theta\FKG\hspace{1.5cm}
                             \mbox{for all}\ F,G\in B(\Omega).
   \end{equation}
   Because $\mu_\theta\circ M_{ij}^{-1}$ and
   $\mu_\theta \circ (S_i^u)^{-1}$ are absolutely continuous
   with respect to $\mu_\theta$, it follows from (\ref{KHF})
   that if $F, \{F_n\}_n$ are uniformly bounded functions
   such that $\int |F_n-F|\mu_\theta(dp)\to_{n\to\infty}0$, then
   $\int|KF_n-KF| \mu_\theta(dp)\to_{n\to\infty} 0$. Thus,
   by standard density arguments we may and shall assume
   $F$ and $G$ to be continuous.

   Define for each $\eps>0$ the truncated intensity measure
   $\nu_\theta^\eps\equiv{\bf 1}_{(\eps,\infty)}\nu_\theta$, and the
   corresponding Poisson measure $\mu_\theta^\eps$, with expectation operator
   $E^\eps_\theta$. Alternatively, if $X$ is distributed in $\Omega$ according
   the $\mu_\theta$, then $\mu_\theta^\eps$
   is the distribution of $T^\eps X:=(X_i1_{X_i>\eps})_i$, that
   is, $\mu^\eps_\theta=\mu_\theta\circ \left(T^\eps\right)^{-1}$.
   Observe that
   $\forall  \delta>0$,
   \[ \mu_\theta(|T^\eps X-X|>\delta)\le \delta^{-1}\,E_\theta|T^\eps  X-X|
          =\delta^{-1}\,E_\theta\,\sum_{p_i<\eps}p_i
          =\delta^{-1}\,\int_0^\eps\,x\,\nu_\theta(dx)
                          \stackrel{\eps\to 0}\longrightarrow 0\ ,     \]
  implying that the measures
   $\mu_\theta^\eps$\ converge weakly to $\mu_\theta$
   as $\eps\to 0$.

   To prove~(\ref{GKF=FKG}) we first write
   \begin{eqnarray} \label{three}
    |E_\theta\GKF-E_\theta\FKG|
      &\le&|E^\eps_\theta\GKF-E_\theta\GKF|+|E^\eps_\theta\GKF-E^\eps_\theta\FKG|
                                                                        \nonumber\\
      &&\hspace*{4.23cm}+\,|E^\eps_\theta\FKG-E_\theta\FKG|
   \end{eqnarray}
   and conclude that the first and third terms
   in~(\ref{three}) converge to $0$ as $\eps\to 0$ by virtue of
   the weak convergence of $\mu^\eps_\theta$ to $\mu_\theta$ and
   $K$'s Feller property, established in Lemma~\ref{Feller}.
    It thus remains to be shown that, for all $F,G\in B(\Omega)$ and $\eps>0$,
   \begin{equation} \label{approxrevers}
    \lim_{\eps\to 0}
      \left|E_\theta^\eps\GKF-E_\theta^\eps\FKG\right|=0.
   \end{equation}
  \bigskip

    \noindent
    The truncated intensity\ \ $\nu_\theta^\eps$\ \ has finite mass\ \ \
    $V_\theta^\eps=\theta\int_\eps^\infty\,{x^{-1}e^{-x}}{}\,dx$,\ \ \
    and thus $N(\RR{}_+)<\infty$,\ \ \ $\mu_\theta^\eps$--a.s.
    In particular each $F\in B(\Omega)$ can be naturally represented as a
    sequence $(F_n)_{n=0}^\infty$ of symmetric $F_n$'s\
    $\in B\left({\bf R}_+^n\right)$, with $\|F_n\|_\infty\le \|F\|_\infty$
    for each $n$. As a result, and in terms of the
    expectation operators $E_{\theta,n}^\eps$ of $\mu_\theta^\eps$
    conditioned on $\{N(\RR{}_+)=n\}$, we may write
    \begin{equation} \label{FG-GF}
      E_\theta^\eps\GKF-E_\theta^\eps\FKG
       =e^{-V_\theta^\eps}\sum_{n=1}^\infty \frac{\left(V_\theta^\eps\right)^n}{n!}
           \left[E_{\theta,n}^\eps\GKF^{^{^{^{}}}}
                            -E_{\theta,n}^\eps\FKG\right],
    \end{equation}
    while by the
    definition~(\ref{KH}) of $\KH$ and the properties stated above of
    the Poisson random measure conditioned on $\{N(\RR{}_+)=n\}$,
    \begin{eqnarray}  \label{EKFG}
 \lefteqn{\frac{\left(V_\theta^\eps\right)^n}{n!}\,
                  E_{\theta,n}^\eps\GKF=}   \nonumber\\
         &&\ \ \ \frac{\beta_m\theta^n}{n!}\,\sum_{\stackrel{i,j=1}{i\ne j}}^n
              \int_\eps^\infty\cdots\int_\eps^\infty \hx_i\,\hx_j
                              F_{n-1}(M_{ij}{\bf x})G_n({\bf x})
                      e^{-|{\bf x}|} \frac{dx_1}{x_1}\,\ldots\,\frac{dx_n}{x_n}
                                                                \nonumber\\
          &&+\frac{\beta_s\theta^n}{n!}\,\sum_{i=1}^n
               \int_\eps^\infty\cdots\int_\eps^\infty \hx_i^2
                     \left(\int_0^1 F_{n+1}(S_i^u{\bf x})G_n({\bf x})\,du\right)
                        e^{-|{\bf x}|}\frac{dx_1}{x_1}\,\ldots\,\frac{dx_n}{x_n}
                                                                \nonumber\\
          &&+\frac{\theta^n}{n!}\,\int_\eps^\infty\cdots\int_\eps^\infty
                \left(1-\beta_m+(\beta_m-\beta_s)
                      \sum_{i=1}^n\hx_i^{\,2}\right) F_n({\bf x})G_n({\bf x})
                     \,e^{-|{\bf x}|}\frac{dx_1}{x_1}\,\ldots\,\frac{dx_n}{x_n}
                                                      \nonumber\\ \nonumber\\
          &&\hspace*{4cm}=: I^{(1)}_n(F,G)+I^{(2)}_n(F,G)+I^{(3)}_n(F,G)\ \ ,
     \end{eqnarray}
where ${\bf x}=(x_1,\ldots,x_n)$.
  Our goal is to prove that this expression, after summing in $n$, is
  roughly symmetric in $F$ and $G$ (as stated precisely
  in~(\ref{approxrevers})). Obviously $I^{(3)}_n(F,G)=I^{(3)}_n(G,F)$,
  and in addition  we aim at showing that
  $I_{n-1}^{(2)}(G,F)\approx I_n^{(1)}(F,G)$\ (with an error appropriately
  small as $\eps\to 0$). This will be achieved by a simple change of variables,
  including the splitting coordinate $u$ in  $I^{(2)}$.
 \bigskip

  In the integral of the $i$--th term in $I^{(2)}_{n-1}(G,F)$
  perform the change of variables  $(u,x_1,\ldots,x_{n-1})\to(y_1,\ldots,y_n)$
  given by ${\bf y}=S_i^u{\bf x}$\ \ \ \ \
  (or\ \ \ \ $(u,{\bf x})=(\frac{y_i}{y_0+y_i},M_{i\,n}{\bf y})$).
  More precisely, $\left\{\begin{array}{l}
                 y_i=ux_i                     \\
                 y_j=x_j,\hspace{.3cm}j\ne i  \\
                 y_n=(1-u)x_i          \\ \end{array}\right.$\\ 
  for which \ \ $|{\bf y}|=|{\bf x}|$\ \ \ \ and\ \ \ \ \
  $dy_1\ldots\,dy_n=x_i\,du\,dx_1\ldots\,dx_{n-1}$,\ \ \
  so that\\

 \noindent
 $I^{(2)}_{n-1}(G,F)=$
  \[  \hspace*{1.3cm}=\frac{\beta_s\theta^{n-1}}{(n-1)!}\sum_{i=1}^{n-1}
        \int_\eps^\infty\cdots\int_\eps^\infty G_n({\bf y})
                                                      F_{n-1}(M_{i\,n}{\bf y})
               \frac{e^{-|{\bf y}|}dy_1\,\ldots\,dy_n}{|{\bf y}|^2\ y_1\ldots\,
                         \breve{y_i}\ldots y_{n-1}}
               +C_n^\eps\]
 ($C_n^\eps$\ is as the term preceding it but with the\
           $dy_i$\ and\ $dy_n$\ integrals taken in\ $[0,\eps]$, and the
notation $\breve{y_i}$ means that the variable $y_i$ has been eliminated from
the denominator)
  \begin{equation} \label{I2}
     = \frac{\beta_s\theta^{n-1}}{(n-2)!}
             \int_\eps^\infty\cdots\int_\eps^\infty G_n({\bf y})
                                                      F_{n-1}(M_{1\,n}{\bf y})
               \frac{e^{-|{\bf y}|}dy_1\,\ldots\,dy_n}
                                {|{\bf y}|^2\ y_2\ldots\, y_{n-1}}
        +C_n^\eps
  \end{equation}
  \medskip

  \noindent
  (by $F_{n-1}$'s symmetry, the sum's $(n-1)$ terms are equal,
  hence the last equality).\\
  \medskip

  On the other hand, and for the same reason of symmetry, the
  $n(n-1)$ terms in $I_n^{(1)}(F,G)$ are all equal so that
  \begin{equation} \label{I1}
   I_n^{(1)}(F,G)=\frac{\beta_m\theta^n}{(n-2)!}
           \int_\eps^\infty\ldots\int_\eps^\infty\,
                  F_{n-1}(M_{1\,2}{\bf x})G_n({\bf x})\,
                  \frac{e^{-|{\bf x}|}\
                  dx_1\ldots dx_n}{|{\bf x}|^2\,x_3\ldots x_n}\ .
  \end{equation}
  Comparing~(\ref{I2}) with~(\ref{I1}), and observing that by definition
  $\beta_m\theta=\beta_s$, we conclude that there exists a $C>0$ such that,
  for $n\geq 2$,
  \begin{eqnarray} \label{Cenbound}
   \lefteqn{
    \left|C_n^\eps\right|:=
   \left|I_{n-1}^{(2)}(G,F)-I_n^{(1)}(F,G)\right|
                                      }\nonumber\\
       &\le& \frac{\|F\|_\infty\|G\|_\infty\beta_s\theta}{(n-2)!}
            \int_0^\eps\!\int_0^\eps\,dy_1\,dy_n\,\frac{1}{((n-2)\eps)^2}
           \left(\theta\,\int_\eps^\infty\,\frac{e^{-y}}{y}\,dy\right)^{n-2}
                                                                \nonumber\\
       &\le&C\,\frac{\left(V_\theta^\eps\right)^{n-2}}{(n-1)!}
   \end{eqnarray}
  Applying~(\ref{Cenbound}) via~(\ref{EKFG}) in~(\ref{FG-GF}) twice, once as
  written and once reversing the roles of $F$ and $G$, and noting that
  $I_1^{(1)}(F,G)=I_1^{(1)}(G,F)=0$, we have
   \begin{eqnarray*}
    \lefteqn{\left|E_\theta^\eps\GKF-E_\theta^\eps\FKG\right|}             \\
        &&\le e^{-V_\theta^\eps}\,\left(
             \sum_{n=2}^\infty\left|I_{n-1}^{(2)}(G,F)-I_n^{(1)}(F,G)\right|
            +\sum_{n=2}^\infty\left|I_{n-1}^{(2)}(F,G)-I_n^{(1)}(G,F)\right|
                                                                    \right) \\
        &&\le 2Ce^{-V_\theta^\eps}\sum_{n=2}^\infty
                           \frac{\left(V_\theta^\eps\right)^{n-2}}{(n-1)!}
            \le \frac{2C}{V_\theta^\eps}
   \end{eqnarray*}
   from which~(\ref{approxrevers}) follows immediately
since $V_\theta^\eps\to_{\eps\to 0} \infty$.
\end{proof}
\begin{proof}[Proof of Theorem \ref{theo-uniq}]
 (a) The Poisson-Dirichlet law $\mu=\widehat{\mu}_\theta$ is reversing by Theorem \ref{reversible}, and
hence invariant. We now show that
it belongs to  $\mathcal A$. Note first that
$\mu_k$ is absolutely continuous with respect to $\mbox{leb}^k$:
for any
$D\subset \Omega_<^k$ with $\mbox{leb}^k(D)=0$, it holds that
$$\mu_k(D) \leq \int_{\RR{}_+} \nu_\theta\left[
\exists\ \jj\in \mathbb{N}^k_{\neq}:
  (X_{j_1},\ldots,
X_{j_k})\in x D\right] d\gamma_\theta(x)=0\,,$$
where we used the fact that under $\mu_\theta$,
$\pi(X)=X/|X|$ and
$|X|$ are independent, with $|X|$ being distributed according to the
Gamma law $\gamma_\theta(dx)$ of density
$1_{x\geq 0}x^{\theta-1} e^{-x}/\Gamma(\theta)$ (see \cite{Ki1}).
It thus suffices to compute the limit
$$p_k(x_1,\ldots,x_k):=\lim_{\delta\to 0}\frac{E_{\widehat{\mu}_\theta}
\left[\# \left\{{\jj}\in \mathbb{N}^k_{\neq} :\,
p_{j_i}\in \left(x_i,x_i+\delta\right)
\,, i=1,\ldots,k\right\}\right]}{\delta^k}\,,$$
where all $x_i$ are pairwise distinct and nonzero, to have
$$m_k(x_1,\ldots,x_k)= p_k(x_1,\ldots,x_k)\prod_{i=1}^kx_i\,.$$
For such $x_1,\ldots,x_k$,
set $I_i^\delta=(x_i,x_i+\delta)$ and $I^\delta=\cup_{i=1}^k I_i^\delta$.
Define
$$L_X^\delta:=\sum_i X_i 1_{\{X_i\not\in I^\delta\}}\,,\quad
N_{x_i}^\delta=\# \{j: X_j\in I_i^\delta\}\,.$$
By the memoryless property of the Poisson process, for any Borel subset
 $A\subset\reals$,
\begin{equation}
\label{limp}
\lim_{\delta\to 0}
P(L_X^\delta\in A\,|\, N_{x_i}^\delta, i=1,\ldots,k)=
P(|X|\in A)=\gamma_\theta(A)\,,
\end{equation}
\abel{limp}
where (\ref{limp}),
as above, is  due to \cite{Ki1}.
Further, recall that $N$ and $(\widehat{X}_i)_i$
are independent. Recall that the density of the Poisson process
at $(y_1,\ldots,y_k)$ is $\theta^k e^{-|y|}/\prod_{i=1}^k y_i$,
where $|y|=y_1+\ldots+y_k$.
Performing
the change of variables $y_i/(z+|y|)=x_i$,
one finds that the  Jacobian of this change of coordinate is
$(z+|y|)^k/(1-|x|)$
 (in computing this Jacobian, it is useful to first
make the change of coordinates $(y_1,\ldots,y_{k-1},|y|)
\mapsto (\bar x_1,\ldots,\bar x_{k-1}, |\bar x|)$
where $|y|,|\bar x|$ are considered as independent coordinates, and note the
block-diagonal structure of the Jacobian).
It follows that
\[
m_k(x_1,\ldots,x_k)=\frac{\theta^k}{(1-|x|)}
 \int_0^\infty \exp\left(-z|x|/(1-|x|)\right)
\gamma_\theta(dz)
=\theta^k (1-|x|)^{\theta-1},
\]
which is real analytic on $\{x\in \R^k: |x|<1\}$.
Thus,
$\widehat{\mu}_\theta\in {\mathcal A}$.
In passing, we note that $m_k(\cdot)=1$ on $\bar \Omega_<^k$ when
$\theta=1$.

\noindent
(b) 1)
First we show that  the family of functions $(m_k)_{k\geq 1}$ associated with $\mu$,
determines $\mu$. To this end, define for  $\jj\in\N^k\ (k\in\N)$  functions
$g_{\jj},\hat{g}_{\jj}:\Omega_1\to [0,1]$  by
\[g_{\jj}(p):=
\sum_{\ii\in\N_{\ne}^k} \prod_{\ell=1}^k
p_{i_\ell}^{j_\ell}\quad\mbox{and}\quad
\hat{g}_{\jj}(p):=\prod_{\ell=1}^k Z_{j_\ell}(p)\quad
 \mbox{where}\quad Z_j(p):=\sum_{i}p_i^j.
\]
Note that  any function $\hat{g}_{\jj}$ with
$\jj\in\N^k$
can be written after expansion of  the product as a (finite) linear combination of functions $g_{\hh}$ with
$\hh\in\N^n, n\geq 1$.
Since we have by the definition of $\mu_k$ that
\begin{equation}
\int
g_{\jj}\ d\mu
=\int_{\bar \Omega_<^k}\prod_{\ell=1}^k x_\ell^{j_\ell-1}\ d\mu_k(x)
=\int_{\bar \Omega_<^k} m_k(x)\prod_{\ell=1}^k x_\ell^{j_\ell-1}\ dx\,,
\label{first}
\end{equation}
\abel{first}
 the
 family $(m_k)_{k\geq 1}$ therefore determines the expectations
$\int \hat{g}_{\jj}\ d\mu\ (\jj\in\N^k, k\geq 1)$.
Consequently, $(m_k)_{k\geq 1}$
determines also the joint laws of the random variables
$(Z_1,\ldots,Z_k),\ k\geq 1,$ under $\mu$. We claim that these laws
characterize $\mu$. Indeed, let $\bar \mu$ be the distribution of
the  random variable $\pi:=(Z_n)_{n\geq 0}: \Omega_1\to
[0,1]^\N$ under $\mu$.
Since
$\pi$ is injective it suffices to show that the distributions of
$(Z_1,\ldots,Z_k),\ k\geq 1,$ under $\mu$ determine $\bar \mu$. But,
since any continuous test
function $F$ on the compact space $[0,1]^\N$ can be uniformly
approximated by the
local function $F_k((x_n)_{n\geq 1}):=F(x_1,\ldots,x_k,0,\ldots)$,
this is true due to
\[\int F\ d\bar \mu=\lim_{k\to\infty}\int F_k\ d\bar
\mu=\lim_{k\to\infty}\int F_k(Z_1,\ldots,Z_k,0,\ldots)\ d\mu.\]

\noindent
2) 
For $\mu\in {\mathcal A}$, the set of numbers
\begin{equation}
\label{corec}
m_k^{({\bfn})}:= m_k^{({\bfn})}(x_1,\ldots,x_k)\Big|_{0,0,\ldots,0}:=
\frac{\partial^{\bf n} m_k}{\partial x_1^{n_1}\cdots \partial x_k^{n_k}}
\Big|_{0,0,\ldots,0}
\end{equation}
\abel{corec}
with $k\geq 1$ and $n_1\geq n_2\geq\ldots\geq n_k\geq 0$
yare enough to characterize $(m_k)_k$, and hence by the first part of the
proof
of b), to characterize $\mu$. It is thus enough
to prove that $K_{\beta_m,\beta_s}$ uniquely determines these numbers.
Toward this end, first note that
\begin{equation}\label{heat}
\int_0^1 m_1(x)\ dx=\mu_1[[0,1]]=1.
\end{equation}\abel{heat}
To simplify notations, we define $m_0\equiv 1$ and extend $m_k$ to a function on $[0,1]^k$ by setting
it 0 on the complement of $\bar \Omega_<^k$. For $k\geq 1$ we have
\begin{equation}
\label{fifth}
\int_0^1
m_k(x_1,\ldots,x_k)dx_1=\left(1-\sum_{i=2}^k x_i\right)m_{k-1}(x_2,\ldots,x_k).
\end{equation}\abel{fifth}
Indeed, for $k=1$ this is (\ref{heat}) while for $k\geq 2$,
and arbitrary $B\in{\mathcal B}_{\bar{\Omega}_<^{k-1}}$,
\begin{eqnarray*}\lefteqn{
\int_B\int_0^1 m_k(x_1,x_2,\ldots,x_k)\ dx_1\ dx_2\ldots dx_k= \mu_k[[0,1]\times B]}\\
&=&E_\mu \Bigg[\sum_{(j_2,\ldots,j_k)\in\mathbb{N}^{k-1}_{\neq}}
     \left(\prod_{i=2}^kp_{j_i}\right)\,{\mathbf
1}_{B}(p_{j_2},\ldots,p_{j_k})\sum_{j_1\notin
\{j_2,\ldots,j_k\}}p_{j_1}1_{[0,1]}(p_{j_1})\Bigg]\\
&=& E_{\mu_{k-1}}\left[1_B(p_2,\ldots,p_k)\left(1-\sum_{i=2}^k p_i\right)\right]\\
& =&
\int_B \left(1-\sum_{i=2}^k x_i\right)m_{k-1}(x_2,\ldots,x_k)\ dx_2\ldots dx_k,
\end{eqnarray*}
which implies (\ref{fifth}).
Now we fix $k\geq 1$, apply
$K_{\beta_m,\beta_s}$ to the test function
$\#\{\jj\in \mathbb{N}^k_{\neq}:
p_{j_i}\in (x_i,x_i+\delta)\,,i=1,\ldots,k\} \delta^{-k}$,
with $(x_1,\ldots,x_k)\in \bar\Omega^k_<$ having pairwise distinct coordinates,
and take $\delta\searrow 0$, which yields the basic relation
\begin{eqnarray*}
&&\beta_m \sum_{i=1}^k
\int_0^{x_i} z(x_i-z)p_{k+1}(x_1,\ldots,x_{i-1},z,x_i-z,x_{i+1},\ldots,x_k)
dz\\
&&+\beta_s \sum_{i=1}^k
\int_{x_i}^1 zp_{k}(x_1,\ldots,x_{i-1},z,x_{i+1},\ldots,x_k)
dz\\
&=&
\beta_m \left(\sum_{i=1}^k x_i(1-x_i)\right) p_k(x_1,\ldots,x_k)+
\beta_s(\sum_{i=1}^k  x_i^2) p_k(x_1,\ldots,x_k)\,.
\end{eqnarray*}
Here the left hand side represents mergings and splittings that produce
a new part roughly at one of the $x_i$-s; the
right hand side represents parts near one of the $x_i$-s
that merge or split.
After multiplying by $x_1\cdots x_k$, rearranging and using (\ref{fifth})
to get rid
of the integral with upper limit 1, we obtain
the equality
\begin{eqnarray}
&&\beta_m \sum_{i=1}^k
x_i\int_0^{x_i} m_{k+1}(x_1,\ldots,x_{i-1},z,x_i-z,x_{i+1},\ldots,x_k)
dz\label{frogs} \\
&&-\beta_s \sum_{i=1}^k
x_i \int_0^{x_i} m_{k}(x_1,\ldots,x_{i-1},z,x_{i+1},\ldots,x_k)
dz\label{cows} \\
&&+\beta_s\sum_{i=1}^k x_i
m_{k-1}(x_1,\ldots,x_{i-1},x_{i+1},\ldots,x_k)\label{magrefa}\\
&&-\beta_s \sum_{i=1}^k\sum_{j=1, j\ne i}^k x_i x_j
m_{k-1}(x_1,\ldots,x_{i-1},x_{i+1},\ldots,x_k)\label{united}\\
&=&\label{second}
\beta_m \left(\sum_{i=1}^k x_i\right) m_k(x_1,\ldots,x_k)+
(\beta_s-\beta_m)(\sum_{i=1}^k  x_i^2) m_k(x_1,\ldots,x_k)\,.
\end{eqnarray}
\abel{frogs}\abel{cows}\abel{magrefa}\abel{united}\abel{second}
We now proceed to show how (\ref{frogs}) -- (\ref{second}) yield all the required
information. As starting point for a recursion,
we show how to compute
$m_k(0,\ldots,0)$ for all $k\geq 1$.
Taking in (\ref{frogs}) -- (\ref{second})
all $x_i\to 0$ except for $x_1$ and using the continuity of the
functions $m_k$ yields
\begin{eqnarray*}
\nonumber
&&\beta_m  \int_0^{x_1} m_{k+1}(z,x_1-z, 0,\ldots,0) dz
-\beta_s \int_0^{x_1} m_k(z,0,\ldots,0) dz\\
&&+\beta_s m_{k-1}(0,\ldots,0) \\
&=& \nonumber
\beta_m  m_k(x_1,0,\ldots,0)+
(\beta_s-\beta_m) x_1 m_k(x_1,0,\ldots,0)\,.
\end{eqnarray*}
Letting $x_1\to 0$  we get $\beta_m m_{k}(0,\ldots,0)=\beta_s m_{k-1}(0,\ldots,0)$.
With $m_0=1$ as start of the recursion this implies
\begin{equation}\label{hit}
m_k(0,\ldots,0)=\theta^k\quad (k\geq 0).
\end{equation}\abel{hit}
For the evaluation of the derivatives of $m_k$
we proceed inductively. Recall the functions $m_k^{({\bf n})}(x_1,\ldots,x_k)$
defined in (\ref{corec}), and write $m_k^{(n_1,n_2,\ldots,n_j)}$,
$j<k$, for $m_k^{(n_1,n_2,\ldots,n_j,0\ldots,0)}$.
Fix ${\bf n}$ such that $n_1\geq n_2\geq \ldots\geq n_k$, with $n_1\geq 2$.
Our analysis rests upon differentiating (\ref{frogs}) --
 (\ref{second}) $n_1$ times with respect to
$x_1$; to make this differentiation easy, call a term
a {\it G term of degree $\ell$} if it is a linear combination of terms
of the form
$$ x_1\int_0^{x_1} m_{k+1}^{(\ell+1)}(z,x_1-z,x_2,\ldots,x_k)dz$$
and
$$ \int_0^{x_1} m_{k+1}^{(\ell)}(z,x_1-z,x_2,\ldots,x_k)dz$$
and
$$ m_{k+1}^{(\ell-1)}(x_1,0,x_2,\ldots,x_k)$$
and
$$ x_1 m_{k+1}^{(\ell)}(x_1,0,x_2,\ldots,x_k)\,.$$
Note that (\ref{frogs}) --
 (\ref{second}) contains one G term of
degree $-1$ in (\ref{frogs}) and that differentiating a G term of degree $\ell$ once yields
a  G term of degree $\ell+1$. Thus, differentiating the G term in (\ref{frogs})
 $n_1\geq 2$ times
and substituting $x_1=0$, we recover a constant multiple
of $m_{k+1}^{(n_1-2)}(0,x_2,\ldots,x_k,0)$.
Similarly, call a term
an {\it H term of degree $\ell$} if it is a linear combination of terms of the form
$$ m_{k}^{(\ell)}(x_1,\ldots,x_k)\quad
\mbox{and}\quad
 x_1 m_{k}^{(\ell+1)}(x_1,\ldots,x_k)\quad
\mbox{and}\quad
 x_1^2 m_{k}^{(\ell+2)}(x_1,\ldots,x_k).$$
Observe, that differentiating an H term of degree $\ell$ produces an
H term of degree $\ell+1$.  If we differentiate twice the term
$x_1\int_0^{x_1} m_k(z,x_2,\ldots,x_k) dz$ in (\ref{cows}) we get an
H term of degree 0. Therefore differentiating this term $n_1\geq 2$ times
results in an H
term of degree $n_1-2$. Since the term  $x_1^2 m_k(x_1,\ldots,x_k)$ in
(\ref{second}) is an H term of degree -2, differentiating this term $n_1$ times
produces also an H term of degree $n_1-2$. Thus both terms produce after $n_1$-fold
differentiation and evaluation at $x_1=0$ a constant multiple of
$m_{k}^{(n_1-2)}(0,x_2,\ldots,x_k)$.
The H term $x_1 m_k(x_1,\ldots,x_k)$  in (\ref{second}) is treated more carefully.
It is easy to see by induction that its $n_1$-th derivative equals $n_1
m_k^{(n_1-1)}(x_1,\ldots,x_k)+ x_1 m_k^{(n_1)}(x_1,\ldots,x_k)$.
Evaluating it at $x_1=0$ gives $n_1
m_k^{(n_1-1)}(0,x_2,\ldots,x_k)$.

Moreover, the terms in (\ref{magrefa}) and (\ref{united}) for $i=1$ vanish when
differentiated twice with respect to $x_1$.

Summarizing the above, we conclude
by differentiating
 (\ref{frogs}) -- (\ref{second}) $n_1\geq 2$ times with respect to $x_1$ and
subsequent
evaluation at
$x_1=0$ that there are some constants
$C_i(n_1)$, such that
\vfill
\eject
\begin{eqnarray*}\nonumber
\lefteqn{\beta_m n_1 m_k^{(n_1-1)}(0,x_2,\ldots,x_k)
\qquad[(\ref{second}a),i=1]}\\
&=&\nonumber
 C_1 m_{k+1}^{(n_1-2)}(0,x_2,\ldots,x_k,0)\qquad[(\ref{frogs}),i=1]\\
&& +
 C_2 m_k^{(n_1-2)}(0,x_2,\ldots,x_k)\qquad[(\ref{second}b),i=1 + (\ref{cows}),i=1]\\
&&\label{lauf}
-\left[\beta_m \left(\sum_{i=2}^k x_i\right)+(\beta_m-\beta_s)
\left(\sum_{i=2}^k x_i^2\right)\right]
m_k^{(n_1)}(0,x_2,\ldots,x_k)\qquad[(\ref{second})]\\
&& +\nonumber
\beta_m \sum_{i=2}^k x_i \int_0^{x_i}
m_{k+1}^{(n_1)}(0,x_2,\ldots,x_{i-1},z,x_i-z,x_{i+1},\ldots,x_k)dz
\qquad[(\ref{frogs})]\\
&& -
\beta_s \sum_{i=2}^k x_i \int_0^{x_i}
m_{k}^{(n_1)}(0,x_2,\ldots,x_{i-1},z,x_{i+1},\ldots,x_k)dz
\qquad[(\ref{cows})]\nonumber
\\
&&+\beta_s\sum_{i=2}^k x_i\nonumber
m_{k-1}^{(n_1)}(0,x_2,\ldots,x_{i-1},x_{i+1},\ldots,x_k)\qquad[(\ref{magrefa})]\\
&&-\beta_s n_1 \sum_{i=2}^k x_i
m_{k-1}^{(n_1-1)}(0,x_2,\ldots,x_{i-1},x_{i+1},\ldots,x_k)
\qquad[(\ref{united}), j=1]\nonumber\\
&&-\beta_s \sum_{i=2}^k\sum_{j=2, j\ne i}^k x_i x_j
m_{k-1}^{(n_1)}(x_1,\ldots,x_{i-1},x_{i+1},\ldots,x_k)\qquad[(\ref{united})]
\end{eqnarray*}\abel{lauf}
For $x_2=\ldots=x_k=0$ only the first three lines do not vanish
and give a recursion which allows us to compute starting with
(\ref{hit}) all
derivatives $m_k^{(n)}(0,\ldots,0)\ (n\geq 0)$.

Further differentiating with respect to $x_2,\ldots,x_k$,
one concludes that
there
exist constants $D^i_{\bfn,\bfn'}$ such that
\begin{eqnarray}
\label{shweishwei}
\beta_m n_1  m_k^{(n_1-1,n_2,\ldots,n_k)}&=&
\!\!\!\!\!\!\!
\sum_{\bfn':|\bfn'|\leq |\bfn|-2, n_i'\leq n_i}
[D^1_{\bfn,\bfn'} m_k^{(\bfn')}+D^2_{\bfn,\bfn'} m_{k+1}^{(\bfn',0)}
+D^3_{\bfn,\bfn'} m_{k-1}^{(\bfn')}]
\nonumber \\
&+&
\!\!\!\!\!\!\!
\!\!\!\!\!\!\!
\!\!\!\!\!\!\!
\sum_{\bfn':|\bfn'|\leq |\bfn|-1, n_i'\leq n_i,n_1=n_1'}
[D^4_{\bfn,\bfn'} m_k^{(\bfn')}+D^5_{\bfn,\bfn'} m_{k+1}^{(\bfn',0)}
+D^6_{\bfn,\bfn'} m_{k-1}^{(\bfn')}]\,,
\end{eqnarray}
\abel{shweishwei}
where $D^3,D^6=0$ unless $\bfn'$ possesses at least one component which is zero.
We now compute iteratively any
of the $m_k^{({\bf n})}$, with $n_1\geq n_2\geq \ldots\geq n_k$: first,
substitute in (\ref{shweishwei})
$n_1=n+1,n_2=1$ to compute
$m_k^{(n,1)}$, for all $n,k$. Then, substitute $n_1=n+1,n_2=j$ ($j\leq n$)
to compute
iteratively $m_k^{(n,j)}$ from the knowledge of the
family $(m_k^{(\ell,j')})_{k,\ell,j'<j}$, etc.
More generally,
having computed the terms $(m_k^{(n_1,n_2,\ldots,n_j)})_{j\leq j_0<k}$,
we compute first $m_k^{(n_1,\ldots,n_{j_0},1)}$
by substituting in (\ref{shweishwei}) $\bfn=(n_1+1,n_2,\ldots,n_{j_0},1)$,
and then proceed inductively as above.
\end{proof}

\section{Concluding remarks}
\label{concluding}

\noindent
1)  We of course conjecture the
\begin{conjecture}
\label{conj1}
Part b) of Theorem \ref{theo-uniq} continues to hold true without the assumption
$\mu\in {\mathcal A}$.
\end{conjecture}
It is tempting to use the technique leading to
(\ref{half}) in order to prove the conjecture
by characterizing the expectations with respect to $\mu$ of suitable
test functions.
One possible way to do that is to consider a family of
polynomials defined as follows. Let ${\bf n}=(n_2,n_3,\ldots,n_d)$ be a
finite sequence of nonnegative integers, with $n_d\geq 1$.
We set $|{\bf n}|=\sum_{j=2}^d j n_j$,
i.e. we consider $\bf n$ as representing a partition of $|{\bf n}|$
having $n_j$ parts of size $j$, and no parts of size $1$.
Recall next $Z_j=Z_j(p)=\sum_i p_i^j$ and the {\bf n}-{\it polynomial}
$$P_{{\bf n}}(p)=\prod_{j=2}^d Z_j^{n_j}:
\Omega_{1}\to \reals\,.$$
$|{\bf n}|$ is the {\it degree} of $P_{\bf n}$, and, with ${\bf n}$ and $d$
as above, $d$ is  the {\it maximal monomial degree} of $P_\bfn$.
Because we do not allow partitions with parts of size $1$, it
holds that $P_\bfn\neq P_{\bfn'}$ if $\bfn\neq \bfn'$ (i.e, there exists
a point $p\in \Omega_1$ such that $P_\bfn(p)\neq P_{\bfn'}(p)$).
It is easy to check that
the family of polynomials $\{P_{\bfn}\}$ is separating for
${\mathcal M}_1(\Omega)$.
Letting $\Delta_\bfn$ denote the expected increment (conditioned on $p$)
of $P_\bfn$ after one step of the process, we have that
$\Delta_\bfn$ is uniformly bounded. Hence, by invariance of $\mu$,
$\int \Delta_\bfn d\mu=0$. Expanding this equality, we get that
\begin{eqnarray}
\label{aroch}
&&
\!\!
\!\!\!\!\!\!
\frac{\beta_m}{\beta_s}
E_\mu\left[ \sum_{\alpha,\beta} p_\alpha p_\beta
\sum_{k=2}^d
\left(\prod_{j=2}^{k-1}(Z^q_{j,\alpha,\beta})^{n_j}\right)
\left(\sum_{\ell=0}^{n_k-1} (Z_k)^\ell
\left(\begin{array}{l}
\!\!n_k\\
\!\!\ell
\end{array}\!\!\right) q_{\alpha,\beta,k}^{n_k-\ell}\right)
\left(\prod_{j=k+1}^d Z_j\right)\right]=
\nonumber\\
&&
\!\!\!\!\!
\!\!\!\!\!
\!\!\!\!\!
-E_\mu\left[
\sum_\alpha p_\alpha^2  \sum_{k=2}^d
\int\!\!\!\left(\!\left(\prod_{j=2}^{k-1}
\!\left(Z^f_{j,\alpha,x}\right)^{n_j}\right)
\!\!\left(\sum_{\ell=0}^{n_k-1} Z_k^\ell\! \left(\begin{array}{l}
\!\!n_k\\
\!\!\ell
\!\!\end{array}\!\!\right)\! f_{\alpha,k,x}^{n_k-\ell}\right)
\!\left(\prod_{j=k+1}^d Z_j^{n_j}\right)
\!\right)\!d\sigma(x)\right]
\nonumber \\
&&
\!\!\!\!\!\!\!\!\!\!\!
+\frac{\beta_m}{\beta_s}
E_\mu \left[\sum_\alpha p_\alpha^2
\sum_{k=2}^d\left(\prod_{j=2}^{k-1}
\left(Z_j+(2^j-2)p_\alpha^j\right)^{n_j}\right)
\right.
\nonumber\\
&&
\;\;\;\;\;\;\;\;\;\;\;
\;\;\;\;\;\;\;\;\;\;\;
\;\;\;\;\;\;\;\;\;\;\;
\left.\left(\sum_{\ell=0}^{n_k-1}Z_k^\ell \left(\begin{array}{l}
\!\!n_k\\
\!\!\ell
\end{array}\!\!\right) \left(
(2^k-2)p_\alpha^k\right)^{n_k-\ell}\right)
\left(\prod_{j=k+1}^d Z_j^{n_j}\right)\right]
\end{eqnarray}
where
$$q_{\alpha,\beta,j}=(p_\alpha+p_\beta)^j-p_\alpha^j-p_\beta^j\geq 0,
\,f_{\alpha,j,x}=[x^j+(1-x)^j-1]p_\alpha^j\leq 0\,,$$
$$
Z^q_{j,\alpha,\beta}=Z_j+q_{\alpha,\beta,j}\,,\,
Z^q_{j,\alpha,x}=Z_j+f_{\alpha,j,x}.$$
Note that all terms in (\ref{aroch}) are positive. Note also that
the right hand side of (\ref{aroch}) is a polynomial of degree
$|\bfn|+2$, with maximal monomial
degree $d+2$, whereas the  left
hand side is a polynomial of degree at most $|\bfn| +2$ and
maximal monomial degree at most $d$.
Let $\pi(k)$ denote the number
of integer partitions of $k$ which do not have parts of size $1$.
Then, there are $\pi(k)$ distinct polynomials of degree $k$,
whereas (\ref{aroch}) provides at most $\pi(k-2)$ relations
between their expected values (involving possibly the
expected value of lower order polynomials). Since always $\pi(k)>\pi(k-2)$,
it does not seem possible to characterize an invariant probability measure
$\mu\in {\mathcal M}_1(\Omega_1)$ using only these algebraic
relations.

\noindent
2) With a lesser degree of confidence we conjecture
\begin{conjecture}
\label{conj2}
For any
$\si\in{\mathcal M}_1((0,1/2])$ and any $\beta_m, \beta_s \in (0,1]$
there exists exactly one $K_{\sigma,\beta_m,\beta_s}$-invariant probability measure
$\mu\in{\mathcal M}_1(\Omega_1)$.
\end{conjecture}

\noindent
3)  We have not been able to resolve whether
the state $\bar{p}=(1,0,0,\ldots)$ is
transient or null-recurrent for $K_{\si,1,1
}$
with $\si=U(0,1/2]$.

\noindent
4) There is much literature concerning
coagulation-fragmentation
processes. Most of the
recent probabilistic
literature deals with processes which exhibit either pure fragmentation
or pure coagulation.
For an extensive review, see \cite{aldous}, and a sample
of more recent references is
\cite{jpda97sac}, \cite{sznitman} and \cite{jpse96cmc}.
Some recent results on coagulation-fragmentation
processes are contained in \cite{jeon}. However, the starting point for
this and previous studies are the coagulation-fragmentation equations,
and it is not clear how to relate those to our model. The functions $m_k$
introduced in the context of  Theorem \ref{theo-uniq} are
related to these equations.

\noindent
5) A characterization of the Poisson-Dirichlet process as the unique measure
coming from an i.i.d. residual allocation model which is
invariant under a {\it split and merge} transformation is given
in \cite{gnedin}. J. Pitman has
pointed out to us that a slight modification
of
this transformation,
preceded
by a size biased permutation and followed by ranking,
is equivalent to our
Markov transition  $K_{\sigma,\beta_m,\beta_s}$. Pitman \cite{pitman}
then
used  this observation
to give an alternative proof of part (a) of Theorem \ref{theo-uniq}.

\noindent
6) Yet another proof of part a) of Theorem \ref{theo-uniq} which avoids the
Poisson representation and Theorem \ref{reversible} can be obtained by
computing the expectation of the polynomials $P_{\bf n}(p)$, defined
in remark 1) above,
under the Poisson-Dirichlet law. We prefer the current proof as it yields more
information and is more transparent.

\noindent
7) A natural extension of Poisson-Dirichlet measures are the {\it two parameter}
Poisson-Dirichlet measures, see e.g. \cite{pityor}.  Pitman raised
the question, which we have not addressed, of whether
there are splitting measures $\sigma$ which would lead to
invariant measures from this family.

\noindent
8) While according to Theorem \ref{theo-uniq}
there is a reversing probability measure
for $\si=U(0,1/2]$
this does not hold for general $\si\in{\mathcal M}_1((0,1/2])$.
For instance, let us assume that
the support of $\si$ is finite.
Then there exist $0<a<b\leq 1/2$ such that $\si[(a,b)]=0$.
To
show that any invariant measure $\mu$ is not reversing it suffices to
find $s,t\in\Omega_1$ such
that the detailed balance equation
\begin{equation}\label{son}
\mu[\{s\}] K_{\si,\beta_m,\beta_s}(s,\{t\})=
\mu[\{t\}] K_{\si,\beta_m,\beta_s}(t,\{s\})
\end{equation}\abel{son}
fails. Due to
Theorem \ref{bit}, $\mu[\{\bar{p}\}]>0$. Now we first refine the partition $\bar{p}$
 by successive splits
until we reach a state
 $p\in\Omega_1$ with $p_1<\eps$, where $\eps>0$ is a small number.
Since $\mu$ has finite support, $\mu[\{p\}]>0$. Then we create from $p$ by
successive mergings some $s\in\Omega_1$ with $a<s_2/s_1<b$, which is possible
if  $\eps$ was chosen
small enough. Again, $\mu[\{s\}]>0$. If we call now $t$ the state which one gets
from $s$ by merging $s_1$ and $s_2$, then the left hand side of (\ref{son}) is
positive. On the other hand, the right hand side of (\ref{son}) is zero because of
$K(t,\{s\})=0$ due to the choice of $a$ and $b$.

\noindent
{\bf Acknowledgment} We thank R. Brooks for suggesting to us the study
of  the fragmentation-coagulation process described here, A.\ Vershik for
the reference to \cite{Tsilevich} and some useful discussions concerning
the Poisson-Dirichlet measure, and A.-S.\ Sznitman for discussions
concerning the uniqueness issue, and for his help with the proof of
Theorem \ref{reversible}. Finally, we thank J. Pitman for
pointing out \cite{gnedin},  making his observation
\cite{pitman} available to us, and for his comments on the literature.

\vspace{1cm}

\begin{tabular}{ll}   Eddy Mayer Wolf & Ofer Zeitouni\\
   Dept.\ of
Mathematics &  Dept.\ of Electrical Engineering\\
Technion, Haifa 32000, Israel&
Technion, Haifa 32000, Israel\\
emw@tx.technion.ac.il&zeitouni@ee.technion.ac.il\\
&www-ee.technion.ac.il/$\!{}_{\mbox{\~{}}}$zeitouni\\
 &\\
 & \\
  Martin Zerner &\\
   Dept.\ of Electrical Engineering&\\
 Technion, Haifa 32000, Israel&\\
  zerner@ee.technion.ac.il&\\
www-ee.technion.ac.il/$\!{}_{\mbox{\~{}}}$zerner&
\end{tabular}

\end{document}